\documentclass[12pt,a4paper]{article}
\usepackage{amsmath,amssymb,amsfonts}
\def\Bbb{\mathbb}

\title{\bf  Dedekind sums take each value infinitely many times}

\author{Kurt Girstmair}
\date{}
\makeatletter
\let\@@maketitle=\maketitle
\def\maketitle{\def\thispagestyle##1{\relax}\@@maketitle}
\makeatother
\newtheorem{theorem}{Theorem}

\newtheorem{lemma}{Lemma}

\def\BE{\begin{equation}}
\def\EE{\end{equation}}
\def\BD{\begin{displaymath}}
\def\ED{\end{displaymath}}
\def\BA{\begin{array}}
\def\EA{\end{array}}
\def\BEA{\begin{eqnarray*}}
\def\EEA{\end{eqnarray*}}
\def\BI{\bibitem}

\def\N{\Bbb N}
\def\Z{\Bbb Z}
\def\Q{\Bbb Q}
\def\R{\Bbb R}

\def\phi{\varphi}

\def\MB{\mbox}
\def\LD{\ldots}
\def\OV{\overline}

\def\MN{\medskip\noindent}
\def\STOP{\hfill$\Box$}

\def\DED{Dedekind }

\begin{document}
\maketitle

\begin{abstract}

\noindent
For $a\in \Z$ and $b\in\N$, $(a,b)=1$, let $s(a,b)$ denote the classical Dedekind sum.
We show that Dedekind sums take this value infinitely many times in the following sense.
There are pairs $(a_i,b_i)$, $i\in\N$, with $b_i$ tending to infinity as $i$ grows, such that $s(a_i,b_i)=s(a,b)$ for all $i\in \N$.

\end{abstract}

%%%%%%%%%%%%%%%%%%%%%%%%%%%%%%%%%%%%%%
\section*{1. Introduction and result}
%%%%%%%%%%%%%%%%%%%%%%%%%%%%%%%%%%%%%%

Let $a$ be an integer, $b$ a natural number, and $(a,b)=1$. The classical \DED sum $s(a,b)$ is defined by
\BD
   s(a,b)=\sum_{k=1}^{b} ((k/b))((ak/b)).
\ED
Here
\BD
  ((x))=\begin{cases}
                 x-\lfloor x\rfloor-1/2 & \MB{ if } x\in\R\smallsetminus \Z; \\
                 0 & \MB{ if } x\in \Z
        \end{cases}
\ED
(see \cite[p. 1]{RaGr}).

Originally, \DED sums appeared in the theory of modular forms (see \cite{Ap2}). But these sums have also interesting applications in
in a number of other fields, so in connection with class numbers, lattice point problems, topology, and algebraic geometry (see  \cite{{At}, {Me}, {RaGr}, {Ur}}).
Starting with Rademacher \cite{Ra}, several authors have studied the distribution of \DED sums (for instance, \cite{{Br}, {Hi}, {Va}}).

It is often more convenient to work with
\BD
 S(a,b)=12s(a,b)
\ED instead. We call $S(a,b)$ a {\em normalized} \DED sum.

Let $r$ be a rational number such that there exist $a\in\Z$ and $b\in\N$, $(a,b)=1$ with $S(a,b)=r$.
Then
\BE
\label{2}
S(a+jb,b)=r\MB{ for all } j\in\Z.
\EE
Accordingly, the value $r$ is taken infinitely many times in a {\em trivial} sense
($b$ fixed, $a$ running through a congruence class mod $b$).

In the present paper, however, we say that the value $r$ is taken {\em infinitely many times} if,
and only if, there exists a sequence $(a_i,b_i)$, $i\in\N$, such that
$b_i\to \infty$ as $i\to\infty$ and $S(a_i,b_i)=r$ for all $i\in\N$.

The only possible $r\in\Z$ that can be the value of a normalized \DED sum is $r=0$. The value $0$ is taken infinitely many times since $S(a,a^2+1)=0$ for all $a\in\N$.
This is well-known (see \cite[p. 28]{RaGr}). Our main result is the following theorem.

\begin{theorem} % Theorem 1 %%%%%%%%%%%%%%%%%%%%%%%
\label{t1}
Let $r\in\Q$ be the value of a normalized \DED sum.  Then the
value $r$ is taken infinitely many times.
\end{theorem} %%%%%%%%%%%%%%%%%%%%%%%%%%%%%%%%%

%%%%%%%%%%%%%%%%%%%%%%%%%%%%%%%%%%%%%%
\section*{2. The proof}
%%%%%%%%%%%%%%%%%%%%%%%%%%%%%%%%%%%%%%

Let $x$ be a real number, $0<x<1$. We consider the regular continued fraction expansion
\BD
  x=[0,c_1,c_2,\LD],
\ED
where the $c_i$ are natural numbers. This expansion is finite, if, and only if, $x\in \Q$. In this case it has the form
\BD
  x=[0, c_1,\LD,c_n]
\ED
with $n\ge 1$ and $c_n\ge 2$. In the present setting the only irrational numbers of interest are quadratic irrationals.
A number $x$ is a quadratic irrational if, and only if, its continued fraction expansion is infinite and periodic. We need only quadratic irrationals that are nearly
purely periodic, i.e.,
\BD
 x=[0,\OV{c_1,\LD,c_L}]=[0,c_1,\LD, c_L,c_1,\LD,c_L,\LD].
\ED
Let $p_k/q_k$, $k\ge 0$, be the $k$th convergent of $x$. The convergents  are defined recursively in a well-known way (see \cite[p. 250]{Hu}).
In particular, $p_k\in\Z$, $q_k\in\N$, and $(p_k,q_k)=1$ for all $k\ge 0$.
Hence  $S(p_k,q_k)$ is the value of a normalized \DED sum.
If $x=a/b=[0,c_1,\LD,c_n]$ is rational, $a, b\in\N$, $(a,b)=1$, then $p_n=a$, $q_n=b$ (and so $p_n/q_n=x$). Otherwise, $q_k$ tends to infinity for $k\to\infty$.

The core of the proof of Theorem \ref{t1} is the following lemma, which is one of the main results of \cite{Gi2}.

\begin{lemma} % Lemma 1 %%%%%%%%%%%%%%%%%%%%%%%
\label{l1}
Let $x=[0,\OV{c_1,\LD,c_L}]$ be a quadratic irrational with odd period length $L$.
If $k\ge 0$, $k\equiv L-1\mod 2L$, then
\BD
  S(p_k,q_k)=S(p_{L-1},q_{L-1}).
\ED
\end{lemma} %%%%%%%%%%%%%%%%%%%%%%%%%%%%%%%%%

\noindent{Remark.} The constant value $S(p_k,q_k)$, $k\equiv L-1\mod 2L$, takes the form
\BE
\label{4}
  S(p_k,q_k)=\sum_{j=1}^L (-1)^ {j-1}c_j+x+x',
\EE
where $x'$ is the algebraic conjugate of $x$ (see \cite{Gi2}).

\MN
{\em Proof of Theorem \ref{t1}.}
Let $S(a,b)$ be the value of a normalized \DED sum.
By (\ref{2}), we may assume that $0\le a<b$, so $0\le a/b<1$. Due to the remark that precedes Theorem \ref{t1},
we suppose that $S(a,b)\ne 0$. Then $a/b\ne 0$.
Let
\BD
  a/b=[0,c_1,\LD,c_n]
\ED
be the continued fraction expansion of $a/b$. We have $n\ge 1$ and, in particular, $c_n\ge 2$.

{\em Case 1:} $n$ is even. Choose an arbitrary natural number $c$ and define
\BD
  x=[0,\OV{c_1,\LD,c_n,c}].
\ED
So this quadratic irrational has the odd number $L=n+1$ as a period length. By Lemma \ref{l1},
the convergents of $x$ satisfy
\BD
  S(p_k,q_k)=S(p_{L-1},q_{L-1})=S(p_n,q_n) =S(a,b)
\ED
for $k\ge 0$, $k\equiv L-1\mod 2L$.
Since $q_k\to\infty$ for $k\to\infty$, the value $S(a,b)$ is taken infinitely many times.

{\em Case 2:} $n$ is odd. We write
$a/b=[0,c_1,\LD,c_{n-1}, c_n-1, 1]$ and put
\BD
  x=[0,\OV{c_1,\LD,c_{n-1},c_n-1,1,1}].
\ED
So the odd number $L=n+2$ is a period length of $x$. Observe that
$p_{n+1}=a,q_{n+1}=b$. We obtain
\BD
  S(p_k,q_k)=S(p_{L-1},q_{L-1})=S(p_{n+1},q_{n+1}) =S(a,b)
\ED
for $k\ge 0$, $k\equiv L-1\mod 2L$. This gives the same result as in Case 1.
\STOP

\MN
{\em Example.} Let $a=5$, $b=14$. Then $S(a,b)=18/7$ and $a/b=[0,2,1,4]$. Here $n=3$, so Case 2 of the proof applies. Accordingly, we define
$x=[0,\OV{2,1,3,1,1}]=-5/7+\sqrt{226}/14$ and have $L=5$. We obtain $p_4=a, q_4=b$ and, for instance, $p_{14}=4535, q_{14}=12614$,
$p_{24}=4090565$, $q_{24}=11377814$, $p_{34}=3689685095$, $q_{34}=10262775614$, where $14$, $24$ and $34$ are $\equiv L-1\mod 2L$.
Indeed, $S(p_{34},q_{34})=S(p_{24},q_{24})=S(p_{14},q_{14})=S(a,b)=18/7$. From (\ref{4}) we also obtain
$S(a,b)=2-1+3-1+1-10/7=18/7$.

\MN
{\em Remarks.}
1. The proof of Theorem \ref{t1} has exhibited a sequence $(a_i,b_i)$ such that $b_i\to\infty$ and $S(a_i,b_i)=S(a,b)$ for all $i\in\N$. The sequence $b_i$, however,
grows exponentially in $i$. This is a consequence of the exponential growth of the denominators $q_k$ of the convergents of $x$.
Accordingly, the set of the numbers $b_i$ is rather thin within the set $\N$.
In a number of special cases the author could establish a sequence $(a_i, b_i)$ of this kind such that $b_i$ is a polynomial of degree 4 in $i$ --- and so the set of the numbers
$b_i$ is considerably denser in $\N$.

2. It would be interesting to know more about the density of the set of all suitable numbers $b_i$. In our case, the following values of $(a_i,b_i)$ with $b_i<1000$ yield $S(a_i,b_i)=18/7$:
$(5,14),(27,70),(13,119),(31,259),(157,406),(47,455),(293,707),(111,854)$. This suggests that the set of all suitable $b_i$ could be relatively dense in $\N$.

%%%%%%%%%%%%%%%%%%%%%%%%%%%%%%%%%%%%%%%%%%%%%%%%%%%%%
%%%%%%%%%%%%%%%%%%%%%%%%%%%%%%%%%%%%%%%%%%%%%%%%%%%%%%%%%%%%%%%%%%%%%%%%%%

\vspace{0.5cm}
\noindent
Kurt Girstmair            \\
Institut f\"ur Mathematik \\
Universit\"at Innsbruck   \\
Technikerstr. 13/7        \\
A-6020 Innsbruck, Austria \\
Kurt.Girstmair@uibk.ac.at

\end{document}